\magnification=\magstep1

\def\c{\centerline}

\noindent
\def\seq#1{\{#1\}_{n=1}^{\infty}}
\nopagenumbers

\c{\bf On J. Borwein's Concept of Sequentially Reflexive Banach Spaces}
\vskip9pt
\c{by P. \O rno}
\vskip18pt

A Banach space $X$ is reflexive if the Mackey topology $\tau(X^*,X)$
on $X^*$ agrees with the norm topology on $X^*$.  Borwein [B] calls
a Banach space $X$ {\it sequentially reflexive\/} provided that every
$\tau(X^*,X)$ convergent {\it sequence\/} in $X^*$ is norm convergent.
The main result in [B] is that $X$ is sequentially reflexive if every
separable subspace of $X$ has separable dual, and Borwein asks for
a characterization of sequentially reflexive spaces.  Here we answer
that question by proving

\proclaim Theorem. {\sl A Banach space $X$ is sequentially reflexive
if and only if $\ell_1$ is not isomorphic to a subspace of $X$.}

{\bf Proof:} Assume first that $\ell_1$ is not isomorphic to a 
subspace of $X$ and let $\seq{x^*_n}$ be a weak$^*$-null
sequence in $X^*$ for which  
the sequence $\seq{<x^*_n,x_n>}$ converges 
to zero for every weakly null
sequence $\seq{x_n}$ in $X$; by the easy Lemma 2.1 in [B] it is
enough to check that such a sequence $\seq{x^*_n}$ must converge in
norm to zero. If not, by passing to a subsequence we can select
a sequence $\seq{x_n}$ in the unit ball of $X$ with $\seq{<x^*_n,x_n>}$
bounded away from zero.  By passing to a further subsequence, we
can assume by Rosenthal's theorem [R], [D, chpt. XI] on Banach spaces
which do not contain isomorphs of $\ell_1$ that $\seq{x_n}$ is
weakly Cauchy. Since $\seq{x^*_n}$ converges weak$^*$ to zero,
by passing to further subsequences and replacing $\seq{x_n}$
with a subsequence of differences ${{x_{2n}-x_{2n-1}}\over 2}$, we can
assume moreover that $\seq{x_n}$ is weakly null. This contradiction 
completes the proof of the first direction.

To go the other way, suppose that $Y$ is a subspace of $X$ which
is isomorphic to $\ell_1$ and let $\seq{e_n}$ be the image 
of the unit vector basis under some isomorphism from $\ell_1$ onto
$Y$.  Define a bounded linear operator from $Y$ into $L_{\infty}[0,1]$
by mapping $e_n$ to the $n$-th Rademacher function $r_n$.  By the injective
property of $L_{\infty}[0,1]$, this operator extends to a bounded
linear operator $T$ from $X$ into $L_{\infty}[0,1]$. Let $r_n^*$ be
the $n$-th Rademacher function in $L_1[0,1]$ considered as a subspace
of $L_{\infty}[0,1]^*$. Thus the sequence $\seq{r_n^*}$, being equivalent
to an orthonormal sequence in a Hilbert space, converges weakly to zero.
Since $L_{\infty}[0,1]$ has the Dunford-Pettis property (cf. [D, p. 113]),
$\seq{r_n^*}$ converges in the Mackey topology to zero, {\it a fortiori}
$\seq{T^*r_n^*}$ converges $\tau(X^*,X)$ to zero.  But 
$<T^*r_n^*,e_n>=<r_n^*,r_n>=1$, so $\seq{T^*r_n^*}$ does not converge
to zero in norm. 

\vskip9pt

\c{\bf References}

\vskip9pt

\item{[B]} J. Borwein, {\sl Asplund spaces are ``sequentially reflexive",\/}
(preprint).

\item{[D]} J. Diestel, {\sl Sequences and series in Banach spaces,\/}
Springer-Verlag Graduate Texts in Mathematics 92, (1984).

\item{[R]} H.~P.~Rosenthal, {\sl A characterization of Banach spaces 
containing $\ell_1$,\/} {\bf Proc. Nat. Acad. Sci. 71,} (1974), 2411--2413.

\bye